\documentclass[11pt]{amsart}
\theoremstyle{plain}
\newtheorem{theorem}{Theorem}

\newtheorem{corollary}[theorem]{Corollary}

\newtheorem*{acknowledgement}{Acknowledgements}

\theoremstyle{remark}

\newtheorem{remarks}{Remarks}

\newtheorem{example}{Example}

\theoremstyle{definition}
\newtheorem{defn}[theorem]{Definition}

\def \R {\mathbf{R}}
\def \Q {\mathbf{Q}}
\def \Z {\mathbf{Z}}

\def \G {\ifmmode\mathcal{G}^3\else$\mathcal{G}^3$\fi}  
\DeclareMathOperator{\CS}{CS}
\DeclareMathOperator{\csinv}{\Delta S}
\begin{document}

\baselineskip.525cm
\title{Isospectrality and $3$-manifold groups}
\thanks{The author was partially supported by NSF Grant 4-50645}
\author[Daniel Ruberman]{Daniel Ruberman}
\address{Department of Mathematics\newline\indent
Brandeis University \newline\indent
Waltham, MA 02254}
\email{\rm{ruberman@binah.cc.brandeis.edu}}
\date\today
\maketitle
\section{Introduction}
In this note, we explain how a well-known construction of isospectral manifolds
leads to an obstruction to a group being the fundamental group of a closed
$3$-dimensional manifold.  The problem of determining, for a given group $G$,
whether there is a closed $3$-manifold $M$ with $\pi_1(M)\cong G$ is readily
seen to be undecidable; let us write $G\in \G$ if there is such a $3$-manifold. 
A standard conjecture (related to Thurston's geometrization program) states that
the only possible finite groups in \G\ are those which act freely and linearly
on
$S^3$, cf.~\cite{thomas:elliptic}.   Algebraic ideas (which apply to the
analogous problem in high-dimensional topology~\cite{davis-milgram:survey}) come
close to proving this conjecture, although there are groups which are not ruled
out by surgery theory but which do not act linearly.   The new invariant we
construct seems to be related to the classical surgery obstructions, but we do
not know the precise relationship.
\begin{acknowledgement}

I would like to thank Robert Guralnick for pointing out the almost-conjugate
subgroups described in Example A, and Jerry Levine and Doug Park for helpful
conversations. 
\end{acknowledgement}

\section{An invariant of a group, and an obstruction}
Let us say that Riemannian manifolds $(M_1,g_1)$ and $(M_2,g_2)$ are
{\sl isospectral on forms} if the Laplace operators on the ($L^2$-completions
of) $\Omega^k(M_i)$ have the same spectrum (counting multiplicities) for all
$k$.  There is a beautiful construction of pairs of isospectral manifolds, due
to Sunada~\cite{sunada}, which starts with the following group-theoretical
concept.
\begin{defn} Let $G$ be a finite group.  Two subgroups $H$ and $K$ are said to
be {\sl almost-conjugate} if any of the following conditions holds:
\begin{enumerate}
\item Every element in $H$ is conjugate to an element of $K$.
\item For each conjugacy class $C$ in $G$, the intersections $C \cap H$ and
$C \cap K$ have the same cardinality.
\item The permutation representation of $G$ on $G/H$ (acting by left
multiplication) is conjugate to the permutation representation of $G$ on $G/K$.
\end{enumerate}
\end{defn}
It is worth noting that the last version could equally be applied to subgroups
$K$ and $H$ of finite index in a possibly infinite group $G$.

Suppose that $M$ is a Riemannian manifold, and that $H,K$ are subgroups of $G$,
and that $\varphi:\pi_1(M) \to G$ is a surjection.  Let
$M_H$ denote the the covering space of $M$ corresponding to the subgroups
$\varphi^{-1}(H)$, with a similar definition for $M_K$.  Sunada
showed~\cite{sunada,berard:transplantation-I} that if $H$ and $K$ are
almost-conjugate, then $M_H$ and $M_K$ are isospectral on forms.  There are many
examples of almost-conjugate subgroups; we will recall two in the last section.

The other ingredients in our obstruction are two related invariants of a
Riemannian $3$-manifold: the $\eta$-invariant~\cite{aps:I} and the Chern-Simons
invariant~\cite{chern-simons:invariant}.  The $\eta$-invariant is a real-valued
Riemannian invariant, which is defined using the spectrum of the Laplacian on
$1$-forms of a Riemannian $3$-manifold.   In particular, two 
$3$-manifolds which are isospectral on forms must have the {\sl same}
$\eta$-invariant.   The Chern-Simons invariant takes its values in $\R/\Z$; the
main property we will need is that the Chern-Simons invariant is multiplicative
under finite covers.   The $\eta$-invariant determines the Chern-Simons
invariant up to an ambiguity of $1/2 \pmod{\Z}$ via the equation~\cite{aps:II}
\begin{equation}\label{cseta1}
3\eta(M) \equiv 2\CS(M) \pmod{\Z}
\end{equation}
Atiyah-Patodi-Singer identify the ambiguity as follows:
\begin{defn}
For any finitely-generated group $\pi$, let $S(\pi)\in \{0,1\}$ denote the
number (modulo $2$) of two-torsion summands in the abelianization $H_1(\pi)$ of
$\pi$.  For a space $X$, let $S(X) = S(\pi_1(X))= S(H_1(X))$. 
\end{defn}
An elementary application of the universal coefficient theorems shows that
$$
S(\pi) = \dim_{\Z_2}H_1(\pi;\Z_2) - \dim_\Q H_1(\pi;\Q)
$$
The precise relationship between Chern-Simons and $\eta$-invariants is then
\begin{equation}\label{cseta2}
\frac{3}{2}\eta(M) - \CS(M) \equiv \frac{1}{2} S(M) \pmod{\Z}
\end{equation}

With all these preliminaries out of the way, we can now define an invariant of
a group $\pi$.  
\begin{defn}  Let $\pi$ be a finitely generated group, and let $\varphi:\pi \to
G$ be a surjection, where $G$ is a finite group with almost-conjugate subgroups
$H$ and $K$.  Define $\csinv(\pi;\varphi) \in \Z_2$ to be
$$
S(\varphi^{-1}(H)) - S(\varphi^{-1}(K))
$$
If  $\csinv(\pi;\varphi) = 0$ for all possible $\varphi$, then we will say that
$\pi$ satisfies the $\CS{}$ condition.
\end{defn}

\begin{theorem}\label{cscond}
Suppose that $\pi = \pi_1(M^3) $ where $M$ is a closed orientable
$3$-manifold.  Then $\pi$ satisfies the $\CS{}$ condition.
\end{theorem}
\begin{proof}
Given a triple $(G,H,K)$ and a surjection $\varphi:\pi_1(M) \to G$, construct
the manifolds $M_H$ and $M_K$.   By definition, the fundamental group of $M_H$
(resp.~$M_K$) is $\varphi^{-1}(H)$ (resp.~$\varphi^{-1}(K)$), so 
$$
\csinv(\pi,\varphi) = S(M_H) -S(M_K)
$$
Since $M_K$ and $M_H$ are isospectral, they have the same $\eta$-invariant.  
On the other hand, $M_K$ and $M_H$ are both finite covering spaces (of the same
degree) of $M$, so they have the same $\CS$ invariants.  It follows from
equation~\eqref{cseta2} that $S(M_H) =S(M_K) $, and hence that 
$
\csinv(\pi,\varphi)= 0.$
\end{proof}

Theorem~\ref{cscond} can be interpreted as saying that a finite group $G$,
having almost conjugate subgroups $H$ and $K$ with $S(H)\neq S(K)$,
cannot act freely on a homotopy sphere.  Since the
$\CS$ condition is homological in nature, it is reasonable to expect that such
a statement should hold with weaker hypotheses.
\begin{theorem}\label{homology}
Let  $G$ be a finite group, with almost conjugate subgroups $H$ and $K$.  If
$G$ acts freely on a $\Z_2$ homology sphere $M^3$, then $S(H) =S(K)$.
\end{theorem}
\begin{proof}
Let $M_H = M/H$ and $M_K= M/K$ as above.
The argument for Theorem~\ref{cscond} applies to show that $S(M_H)
=S(M_K)$, so we need to see that $S(K) = S(M_K)$, with the analogous
statement holding for $H$.  According to the remark after the definition of
$S$, it suffices to prove that 
$$
 H_1(M_K;\Q) \cong H_1(K;\Q)  (= \{0\})
\quad \mathrm{and}\quad
H_1(K;\Z_2) \cong  H_1(M_K;\Z_2)
$$
The first follows by an elementary transfer argument, and the second follows
directly from the spectral sequences for the (regular) covering $M \to M_K$,
with $\Z_2$ coefficients.  
\end{proof}
\section{An example}
To find groups which are excluded from membership in $\G$ by the $\CS$
condition, one can simply look among the known examples of almost-conjugate
subgroups of finite groups until one finds one for which $S(H) \neq
S(K)$.  

\begin{example}
One example of this phenomenon was pointed out
to me by R.~Guralnick.  Quoting from a letter of several years ago, 
\begin{quote}
Let
G=$M_{23}
$, the Mathieu group of degree 23. There are two subgroups $H=2^4.A_7$ (that is
there is a normal subgroup of order $2^4$ and the quotient is $A_7$) and$
K=L_3(4).2$ which induce the same character (this is evident from the Atlas of
finite groups for example).  $H$ is perfect so $S(H)=0$, while
$K/K' = \Z/2$ so $S(K)=1.$
\end{quote}
{}From this example, one concludes that $M_{23}$
is not the fundamental group of any $3$-manifold, although this is can also be
shown by cohomological means.
\end{example}

\begin{example}
The simplest example I know is based on the following observation. Any group of
order $n$ acts on itself by left multiplication, and hence imbeds as a subgroup
of the symmetric group $S_n$.  Let
$H$ and $K$ be groups of order $n$, which have the same number of elements
of order $k$, for any $k|n$.  Then $H$ and $K$ are almost-conjugate subgroups
of $S_n$. To apply this, let $H = \Z_4 \oplus \Z_2 \oplus \Z_2$, which has
$S(H) = 1$.  In the listing of $2$-groups of small
order~\cite{hall:2-groups} can be found a group denoted $ 16 \Gamma_2 c_1$, with
the same number of elements of any order as does $H$.   (See page 39,
of~\cite{hall:2-groups} for this information, and page 16 for a presentation of
the group.  With this presentation, I verified the necessary algebraic facts
using the computer program GAP~\cite{gap93}.)  Letting
$K =  16\Gamma_2c_1$ which has abelianization $\Z_4 \oplus \Z_2$, we deduce that
the symmetric group $S_{16}$ does not satisfy condition $\CS$.
\end{example}
{}From this example, and the theorem, we deduce
\begin{corollary}
The symmetric group $S_n$, for $n\ge 16$, is not the fundamental group of any
closed orientable $3$-manifold.
\end{corollary}
For $n=16$ this is the conclusion of theorem~\ref{cscond}, while for $n >16$
we use the observation that if $\pi \in \G$, then so is any finite index
subgroup of $\pi$.

It has been known since the work of
Milnor~\cite{milnor:sn} that a symmetric group cannot act freely on a homology
sphere, and hence could not be the fundamental group of a $3$-manifold, so that
the corollary is certainly not new.  Proofs and generalizations (along the
lines of Theorem~\ref{homology}) of Milnor's result were given by
R.~Lee~\cite{lee:semi} and J.~Davis~\cite{davis:semi} using the concept of the
semicharacteristic
$\chi_{\frac{1}{2}}$.   In this vein, it seems quite suggestive that for a
$3$-manifold, we have (for any coefficient field $\mathbf{F}$) by definition
$\chi_{\frac{1}{2}}(M;\mathbf{F}) = \dim_\mathbf{F}(H_0(M;\mathbf{F}) +
\dim_\mathbf{F}(H_0(M;\mathbf{F})$.  Hence
$$
S(M) = \chi_{\frac{1}{2}}(M;\Q) - \chi_{\frac{1}{2}}(M;\Z_2)
$$
However, I do not see how to deduce theorems~\ref{cscond} or~\ref{homology}
using more traditional methods of surgery theory.

Finally, the reader may wonder whether the invariant $\csinv$ can be used to
decide the membership in $\G$ of those finite groups which remain for the
moment out the reach of surgery theory.  According to known theorems (see
again~\cite{davis-milgram:survey}), the key case to decide is that of the
`generalized quaternion groups'
$Q(8a,b,c)$.  Unfortunately, our methods do not shed any light on this
question. 
\begin{theorem}
The groups $Q(8a,b,c)$ all satisfy condition $\CS$.
\end{theorem}
The proof of the theorem is a fairly tedious examination of the subgroup
structure of $Q(8a,b,c)$, and is omitted.
\begin{remarks}
One can write down finite presentations of groups which condition $\CS$
prevents from being the fundamental group of a $3$-manifold.   These examples
would be of greater interest if they satisfied other necessary conditions, for
example if they were Poincar\'e duality groups of dimension $3$.  I would be
interested to know of any examples of $3$-dimensional Poincar\'e duality
groups, which are not already known to be $3$-manifold groups.  Finally,
the proof of Theorem~\ref{cscond} suggests the question of whether the
Chern-Simons invariant is in fact a spectral invariant, i.e., is
determined by the spectrum of the Laplacian.  
\end{remarks}


\newcommand{\etalchar}[1]{$^{#1}$}
\providecommand{\bysame}{\leavevmode\hbox to3em{\hrulefill}\thinspace}

\end{document}